\documentclass[12pt]{article}

\usepackage{amssymb,amsbsy,amsthm,amsmath,amscd,epsfig,graphicx,times}

\edef\thinlines{\the\catcode`@ }%
\catcode`@ = 11
\let\@oldatcatcode = \thinlines

\def\smash@@{\relax 
  \ifmmode\def\next{\mathpalette\mathsm@sh}\else\let\next\makesm@sh
  \fi\next}
\def\makesm@sh#1{\setbox\z@\hbox{#1}\finsm@sh}
\def\mathsm@sh#1#2{\setbox\z@\hbox{$\m@th#1{#2}$}\finsm@sh}
\def\finsm@sh{\ht\z@\z@ \dp\z@\z@ \box\z@}

\edef\@oldandcatcode{\the\catcode`& }%
\catcode`& = 11

\def\&whilenoop#1{}%
\def\&whiledim#1\do #2{\ifdim #1\relax#2\&iwhiledim{#1\relax#2}\fi}%
\def\&iwhiledim#1{\ifdim #1\let\&nextwhile=\&iwhiledim
        \else\let\&nextwhile=\&whilenoop\fi\&nextwhile{#1}}%

\newif\if&negarg
\newdimen\&wholewidth
\newdimen\&halfwidth

\font\tenln=line10

\def\thinlines{\let\&linefnt\tenln \let\&circlefnt\tencirc
  \&wholewidth\fontdimen8\tenln \&halfwidth .5\&wholewidth}%
\def\thicklines{\let\&linefnt\tenlnw \let\&circlefnt\tencircw
  \&wholewidth\fontdimen8\tenlnw \&halfwidth .5\&wholewidth}%

\def\drawline(#1,#2)#3{\&xarg #1\relax \&yarg #2\relax \&linelen=#3\relax
  \ifnum\&xarg =0 \&vline \else \ifnum\&yarg =0 \&hline \else \&sline\fi\fi}%

\def\&sline{\leavevmode
  \ifnum\&xarg< 0 \&negargtrue \&xarg -\&xarg \&yyarg -\&yarg
  \else \&negargfalse \&yyarg \&yarg \fi
  \ifnum \&yyarg >0 \&tempcnta\&yyarg \else \&tempcnta -\&yyarg \fi
  \ifnum\&tempcnta>6 \&badlinearg \&yyarg0 \fi
  \ifnum\&xarg>6 \&badlinearg \&xarg1 \fi
  \setbox\&linechar\hbox{\&linefnt\&getlinechar(\&xarg,\&yyarg)}%
  \ifnum \&yyarg >0 \let\&upordown\raise \&clnht\z@
  \else\let\&upordown\lower \&clnht \ht\&linechar\fi
  \&clnwd=\wd\&linechar
  \&whiledim \&clnwd <\&linelen \do {%
    \&upordown\&clnht\copy\&linechar
    \advance\&clnht \ht\&linechar
    \advance\&clnwd \wd\&linechar
  }%
  \advance\&clnht -\ht\&linechar
  \advance\&clnwd -\wd\&linechar
  \&tempdima\&linelen\advance\&tempdima -\&clnwd
  \&tempdimb\&tempdima\advance\&tempdimb -\wd\&linechar
  \hskip\&tempdimb \multiply\&tempdima \@m
  \&tempcnta \&tempdima \&tempdima \wd\&linechar \divide\&tempcnta \&tempdima
  \&tempdima \ht\&linechar \multiply\&tempdima \&tempcnta
  \divide\&tempdima \@m
  \advance\&clnht \&tempdima
  \ifdim \&linelen <\wd\&linechar \hskip \wd\&linechar
  \else\&upordown\&clnht\copy\&linechar\fi}%

\def\&hline{\vrule height \&halfwidth depth \&halfwidth width \&linelen}%

\def\&getlinechar(#1,#2){\&tempcnta#1\relax\multiply\&tempcnta 8
  \advance\&tempcnta -9 \ifnum #2>0 \advance\&tempcnta #2\relax\else
  \advance\&tempcnta -#2\relax\advance\&tempcnta 64 \fi
  \char\&tempcnta}%

\def\drawvector(#1,#2)#3{\&xarg #1\relax \&yarg #2\relax
  \&tempcnta \ifnum\&xarg<0 -\&xarg\else\&xarg\fi
  \ifnum\&tempcnta<5\relax \&linelen=#3\relax
    \ifnum\&xarg =0 \&vvector \else \ifnum\&yarg =0 \&hvector
    \else \&svector\fi\fi\else\&badlinearg\fi}%

\def\&hvector{\ifnum\&xarg<0 \rlap{\&linefnt\&getlarrow(1,0)}\fi \&hline
  \ifnum\&xarg>0 \llap{\&linefnt\&getrarrow(1,0)}\fi}%

\def\&vvector{\ifnum \&yarg <0 \&downvector \else \&upvector \fi}%

\def\&svector{\&sline
  \&tempcnta\&yarg \ifnum\&tempcnta <0 \&tempcnta=-\&tempcnta\fi
  \ifnum\&tempcnta <5
    \if&negarg\ifnum\&yarg>0                   
      \llap{\lower\ht\&linechar\hbox to\&linelen{\&linefnt
        \&getlarrow(\&xarg,\&yyarg)\hss}}\else 
      \llap{\hbox to\&linelen{\&linefnt\&getlarrow(\&xarg,\&yyarg)\hss}}\fi
    \else\ifnum\&yarg>0                        
      \&tempdima\&linelen \multiply\&tempdima\&yarg
      \divide\&tempdima\&xarg \advance\&tempdima-\ht\&linechar
      \raise\&tempdima\llap{\&linefnt\&getrarrow(\&xarg,\&yyarg)}\else
      \&tempdima\&linelen \multiply\&tempdima-\&yarg 
      \divide\&tempdima\&xarg
      \lower\&tempdima\llap{\&linefnt\&getrarrow(\&xarg,\&yyarg)}\fi\fi
  \else\&badlinearg\fi}%

\def\&getlarrow(#1,#2){\ifnum #2 =\z@ \&tempcnta='33\else
\&tempcnta=#1\relax\multiply\&tempcnta \sixt@@n \advance\&tempcnta
-9 \&tempcntb=#2\relax\multiply\&tempcntb \tw@ \ifnum \&tempcntb >0
\advance\&tempcnta \&tempcntb\relax \else\advance\&tempcnta
-\&tempcntb\advance\&tempcnta 64
\fi\fi\char\&tempcnta}%

\def\&getrarrow(#1,#2){\&tempcntb=#2\relax
\ifnum\&tempcntb < 0 \&tempcntb=-\&tempcntb\relax\fi \ifcase
\&tempcntb\relax \&tempcnta='55 \or \ifnum #1<3
\&tempcnta=#1\relax\multiply\&tempcnta 24 \advance\&tempcnta -6
\else \ifnum #1=3 \&tempcnta=49 \else\&tempcnta=58 \fi\fi\or \ifnum
#1<3 \&tempcnta=#1\relax\multiply\&tempcnta 24 \advance\&tempcnta -3
\else \&tempcnta=51\fi\or \&tempcnta=#1\relax\multiply\&tempcnta
\sixt@@n \advance\&tempcnta -\tw@ \else
\&tempcnta=#1\relax\multiply\&tempcnta \sixt@@n \advance\&tempcnta 7
\fi\ifnum #2<0 \advance\&tempcnta 64 \fi
\char\&tempcnta}%

\def\&vline{\ifnum \&yarg <0 \&downline \else \&upline\fi}%

\def\&upline{\hbox to \z@{\hskip -\&halfwidth \vrule width \&wholewidth
   height \&linelen depth \z@\hss}}%

\def\&downline{\hbox to \z@{\hskip -\&halfwidth \vrule width \&wholewidth
   height \z@ depth \&linelen \hss}}%

\def\&upvector{\&upline\setbox\&tempboxa\hbox{\&linefnt\char'66}\raise
     \&linelen \hbox to\z@{\lower \ht\&tempboxa\box\&tempboxa\hss}}%

\def\&downvector{\&downline\lower \&linelen
      \hbox to \z@{\&linefnt\char'77\hss}}%

\def\&badlinearg{\errmessage{Bad \string\arrow\space argument.}}%

\thinlines

\countdef\&xarg     0 \countdef\&yarg     2 \countdef\&yyarg    4
\countdef\&tempcnta 6 \countdef\&tempcntb 8

\dimendef\&linelen  0 \dimendef\&clnwd    2 \dimendef\&clnht    4
\dimendef\&tempdima 6 \dimendef\&tempdimb 8

\chardef\@arrbox    0 \chardef\&linechar  2
\chardef\&tempboxa  2           

\let\lft^%
\let\rt_

\newif\if@pslope 
\def\@findslope(#1,#2){\ifnum#1>0
  \ifnum#2>0 \@pslopetrue \else\@pslopefalse\fi \else
  \ifnum#2>0 \@pslopefalse \else\@pslopetrue\fi\fi}%

\def\generalsmap(#1,#2){\getm@rphposn(#1,#2)\plnmorph\futurelet\next\addm@rph}%

\def\sline(#1,#2){\setbox\@arrbox=\hbox{\drawline(#1,#2){\sarrowlength}}%
  \@findslope(#1,#2)\d@@blearrfalse\generalsmap(#1,#2)}%
\def\arrow(#1,#2){\setbox\@arrbox=\hbox{\drawvector(#1,#2){\sarrowlength}}%
  \@findslope(#1,#2)\d@@blearrfalse\generalsmap(#1,#2)}%

\newif\ifd@@blearr

\def\bisline(#1,#2){\@findslope(#1,#2)%
  \if@pslope \let\@upordown\raise \else \let\@upordown\lower\fi
  \getch@nnel(#1,#2)\setbox\@arrbox=\hbox{\@upordown\@vchannel
    \rlap{\drawline(#1,#2){\sarrowlength}}%
      \hskip\@hchannel\hbox{\drawline(#1,#2){\sarrowlength}}}%
  \d@@blearrtrue\generalsmap(#1,#2)}%
\def\biarrow(#1,#2){\@findslope(#1,#2)%
  \if@pslope \let\@upordown\raise \else \let\@upordown\lower\fi
  \getch@nnel(#1,#2)\setbox\@arrbox=\hbox{\@upordown\@vchannel
    \rlap{\drawvector(#1,#2){\sarrowlength}}%
      \hskip\@hchannel\hbox{\drawvector(#1,#2){\sarrowlength}}}%
  \d@@blearrtrue\generalsmap(#1,#2)}%
\def\adjarrow(#1,#2){\@findslope(#1,#2)%
  \if@pslope \let\@upordown\raise \else \let\@upordown\lower\fi
  \getch@nnel(#1,#2)\setbox\@arrbox=\hbox{\@upordown\@vchannel
    \rlap{\drawvector(#1,#2){\sarrowlength}}%
      \hskip\@hchannel\hbox{\drawvector(-#1,-#2){\sarrowlength}}}%
  \d@@blearrtrue\generalsmap(#1,#2)}%

\newif\ifrtm@rph
\def\@shiftmorph#1{\hbox{\setbox0=\hbox{$\scriptstyle#1$}%
  \setbox1=\hbox{\hskip\@hm@rphshift\raise\@vm@rphshift\copy0}%
  \wd1=\wd0 \ht1=\ht0 \dp1=\dp0 \box1}}%
\def\@hm@rphshift{\ifrtm@rph
  \ifdim\hmorphposnrt=\z@\hmorphposn\else\hmorphposnrt\fi \else
  \ifdim\hmorphposnlft=\z@\hmorphposn\else\hmorphposnlft\fi \fi}%
\def\@vm@rphshift{\ifrtm@rph
  \ifdim\vmorphposnrt=\z@\vmorphposn\else\vmorphposnrt\fi \else
  \ifdim\vmorphposnlft=\z@\vmorphposn\else\vmorphposnlft\fi \fi}%

\def\addm@rph{\ifx\next\lft\let\temp=\lftmorph\else
  \ifx\next\rt\let\temp=\rtmorph\else\let\temp\relax\fi\fi \temp}%

\def\plnmorph{\dimen1\wd\@arrbox \ifdim\dimen1<\z@ \dimen1-\dimen1\fi
  \vcenter{\box\@arrbox}}%
\def\lftmorph\lft#1{\rtm@rphfalse \setbox0=\@shiftmorph{#1}%
  \if@pslope \let\@upordown\raise \else \let\@upordown\lower\fi
  \llap{\@upordown\@vmorphdflt\hbox to\dimen1{\hss 
    \llap{\box0}\hss}\hskip\@hmorphdflt}\futurelet\next\addm@rph}%
\def\rtmorph\rt#1{\rtm@rphtrue \setbox0=\@shiftmorph{#1}%
  \if@pslope \let\@upordown\lower \else \let\@upordown\raise\fi
  \llap{\@upordown\@vmorphdflt\hbox to\dimen1{\hss
    \rlap{\box0}\hss}\hskip-\@hmorphdflt}\futurelet\next\addm@rph}%

\def\getm@rphposn(#1,#2){\ifd@@blearr \dimen@\morphdist \advance\dimen@ by
  .5\channelwidth \@getshift(#1,#2){\@hmorphdflt}{\@vmorphdflt}{\dimen@}\else
  \@getshift(#1,#2){\@hmorphdflt}{\@vmorphdflt}{\morphdist}\fi}%

\def\getch@nnel(#1,#2){\ifdim\hchannel=\z@ \ifdim\vchannel=\z@
    \@getshift(#1,#2){\@hchannel}{\@vchannel}{\channelwidth}%
    \else \@hchannel\hchannel \@vchannel\vchannel \fi
  \else \@hchannel\hchannel \@vchannel\vchannel \fi}%

\def\@getshift(#1,#2)#3#4#5{\dimen@ #5\relax
  \&xarg #1\relax \&yarg #2\relax
  \ifnum\&xarg<0 \&xarg -\&xarg \fi
  \ifnum\&yarg<0 \&yarg -\&yarg \fi
  \ifnum\&xarg<\&yarg \&negargtrue \&yyarg\&xarg \&xarg\&yarg \&yarg\&yyarg\fi
  \ifcase\&xarg \or  
    \ifcase\&yarg    
      \dimen@i \z@ \dimen@ii \dimen@ \or 
      \dimen@i .7071\dimen@ \dimen@ii .7071\dimen@ \fi \or
    \ifcase\&yarg    
      \or 
      \dimen@i .4472\dimen@ \dimen@ii .8944\dimen@ \fi \or
    \ifcase\&yarg    
      \or 
      \dimen@i .3162\dimen@ \dimen@ii .9486\dimen@ \or
      \dimen@i .5547\dimen@ \dimen@ii .8321\dimen@ \fi \or
    \ifcase\&yarg    
      \or 
      \dimen@i .2425\dimen@ \dimen@ii .9701\dimen@ \or\or
      \dimen@i .6\dimen@ \dimen@ii .8\dimen@ \fi \or
    \ifcase\&yarg    
      \or 
      \dimen@i .1961\dimen@ \dimen@ii .9801\dimen@ \or
      \dimen@i .3714\dimen@ \dimen@ii .9284\dimen@ \or
      \dimen@i .5144\dimen@ \dimen@ii .8575\dimen@ \or
      \dimen@i .6247\dimen@ \dimen@ii .7801\dimen@ \fi \or
    \ifcase\&yarg    
      \or 
      \dimen@i .1645\dimen@ \dimen@ii .9864\dimen@ \or\or\or\or
      \dimen@i .6402\dimen@ \dimen@ii .7682\dimen@ \fi \fi
  \if&negarg \&tempdima\dimen@i \dimen@i\dimen@ii \dimen@ii\&tempdima\fi
  #3\dimen@i\relax #4\dimen@ii\relax }%

\catcode`\&=4  

\def\generalhmap{\futurelet\next\@generalhmap}%
\def\@generalhmap{\ifx\next^ \let\temp\generalhm@rph\else
  \ifx\next_ \let\temp\generalhm@rph\else \let\temp\m@kehmap\fi\fi \temp}%
\def\generalhm@rph#1#2{\ifx#1^
    \toks@=\expandafter{\the\toks@#1{\rtm@rphtrue\@shiftmorph{#2}}}\else
    \toks@=\expandafter{\the\toks@#1{\rtm@rphfalse\@shiftmorph{#2}}}\fi
  \generalhmap}%
\def\m@kehmap{\mathrel{\smash@@{\the\toks@}}}%

\def\mapright{\toks@={\mathop{\vcenter{\smash@@{\drawrightarrow}}}\limits}%
  \generalhmap}%
\def\mapleft{\toks@={\mathop{\vcenter{\smash@@{\drawleftarrow}}}\limits}%
  \generalhmap}%
\def\bimapright{\toks@={\mathop{\vcenter{\smash@@{\drawbirightarrow}}}\limits}%
  \generalhmap}%
\def\bimapleft{\toks@={\mathop{\vcenter{\smash@@{\drawbileftarrow}}}\limits}%
  \generalhmap}%
\def\adjmapright{\toks@={\mathop{\vcenter{\smash@@{\drawadjrightarrow}}}\limits}%
  \generalhmap}%
\def\adjmapleft{\toks@={\mathop{\vcenter{\smash@@{\drawadjleftarrow}}}\limits}%
  \generalhmap}%
\def\hline{\toks@={\mathop{\vcenter{\smash@@{\drawhline}}}\limits}%
  \generalhmap}%
\def\bihline{\toks@={\mathop{\vcenter{\smash@@{\drawbihline}}}\limits}%
  \generalhmap}%

\def\drawrightarrow{\hbox{\drawvector(1,0){\harrowlength}}}%
\def\drawleftarrow{\hbox{\drawvector(-1,0){\harrowlength}}}%
\def\drawbirightarrow{\hbox{\raise.5\channelwidth
  \hbox{\drawvector(1,0){\harrowlength}}\lower.5\channelwidth
  \llap{\drawvector(1,0){\harrowlength}}}}%
\def\drawbileftarrow{\hbox{\raise.5\channelwidth
  \hbox{\drawvector(-1,0){\harrowlength}}\lower.5\channelwidth
  \llap{\drawvector(-1,0){\harrowlength}}}}%
\def\drawadjrightarrow{\hbox{\raise.5\channelwidth
  \hbox{\drawvector(-1,0){\harrowlength}}\lower.5\channelwidth
  \llap{\drawvector(1,0){\harrowlength}}}}%
\def\drawadjleftarrow{\hbox{\raise.5\channelwidth
  \hbox{\drawvector(1,0){\harrowlength}}\lower.5\channelwidth
  \llap{\drawvector(-1,0){\harrowlength}}}}%
\def\drawhline{\hbox{\drawline(1,0){\harrowlength}}}%
\def\drawbihline{\hbox{\raise.5\channelwidth
  \hbox{\drawline(1,0){\harrowlength}}\lower.5\channelwidth
  \llap{\drawline(1,0){\harrowlength}}}}%

\def\generalvmap{\futurelet\next\@generalvmap}%
\def\@generalvmap{\ifx\next\lft \let\temp\generalvm@rph\else
  \ifx\next\rt \let\temp\generalvm@rph\else \let\temp\m@kevmap\fi\fi \temp}%
\toksdef\toks@@=1
\def\generalvm@rph#1#2{\ifx#1\rt 
    \toks@=\expandafter{\the\toks@
      \rlap{$\vcenter{\rtm@rphtrue\@shiftmorph{#2}}$}}\else 
    \toks@@={\llap{$\vcenter{\rtm@rphfalse\@shiftmorph{#2}}$}}%
    \toks@=\expandafter\expandafter\expandafter{\expandafter\the\expandafter
      \toks@@ \the\toks@}\fi \generalvmap}%
\def\m@kevmap{\the\toks@}%

\def\mapdown{\toks@={\vcenter{\drawdownarrow}}\generalvmap}%
\def\mapup{\toks@={\vcenter{\drawuparrow}}\generalvmap}%
\def\bimapdown{\toks@={\vcenter{\drawbidownarrow}}\generalvmap}%
\def\bimapup{\toks@={\vcenter{\drawbiuparrow}}\generalvmap}%
\def\adjmapdown{\toks@={\vcenter{\drawadjdownarrow}}\generalvmap}%
\def\adjmapup{\toks@={\vcenter{\drawadjuparrow}}\generalvmap}%
\def\vline{\toks@={\vcenter{\drawvline}}\generalvmap}%
\def\bivline{\toks@={\vcenter{\drawbivline}}\generalvmap}%

\def\drawdownarrow{\hbox to5pt{\hss\drawvector(0,-1){\varrowlength}\hss}}%
\def\drawuparrow{\hbox to5pt{\hss\drawvector(0,1){\varrowlength}\hss}}%
\def\drawbidownarrow{\hbox to5pt{\hss\hbox{\drawvector(0,-1){\varrowlength}}%
  \hskip\channelwidth\hbox{\drawvector(0,-1){\varrowlength}}\hss}}%
\def\drawbiuparrow{\hbox to5pt{\hss\hbox{\drawvector(0,1){\varrowlength}}%
  \hskip\channelwidth\hbox{\drawvector(0,1){\varrowlength}}\hss}}%
\def\drawadjdownarrow{\hbox to5pt{\hss\hbox{\drawvector(0,-1){\varrowlength}}%
  \hskip\channelwidth\lower\varrowlength
  \hbox{\drawvector(0,1){\varrowlength}}\hss}}%
\def\drawadjuparrow{\hbox to5pt{\hss\hbox{\drawvector(0,1){\varrowlength}}%
  \hskip\channelwidth\raise\varrowlength
  \hbox{\drawvector(0,-1){\varrowlength}}\hss}}%
\def\drawvline{\hbox to5pt{\hss\drawline(0,1){\varrowlength}\hss}}%
\def\drawbivline{\hbox to5pt{\hss\hbox{\drawline(0,1){\varrowlength}}%
  \hskip\channelwidth\hbox{\drawline(0,1){\varrowlength}}\hss}}%

\def\commdiag#1{\null\,
  \vcenter{\commdiagbaselines
  \m@th\ialign{\hfil$##$\hfil&&\hfil$\mkern4mu ##$\hfil\crcr
      \mathstrut\crcr\noalign{\kern-\baselineskip}
      #1\crcr\mathstrut\crcr\noalign{\kern-\baselineskip}}}\,}%

\def\commdiagbaselines{\baselineskip15pt \lineskip3pt \lineskiplimit3pt }%
\def\gridcommdiag#1{\null\,
  \vcenter{\offinterlineskip
  \m@th\ialign{&\vbox to\vgrid{\vss
    \hbox to\hgrid{\hss\smash@@{$##$}\hss}}\crcr
      \mathstrut\crcr\noalign{\kern-\vgrid}
      #1\crcr\mathstrut\crcr\noalign{\kern-.5\vgrid}}}\,}%

\newdimen\harrowlength \harrowlength=60pt
\newdimen\varrowlength \varrowlength=.618\harrowlength
\newdimen\sarrowlength \sarrowlength=\harrowlength

\newdimen\hmorphposn \hmorphposn=\z@
\newdimen\vmorphposn \vmorphposn=\z@
\newdimen\morphdist  \morphdist=4pt

\dimendef\@hmorphdflt 0       
\dimendef\@vmorphdflt 2       

\newdimen\hmorphposnrt  \hmorphposnrt=\z@
\newdimen\hmorphposnlft \hmorphposnlft=\z@
\newdimen\vmorphposnrt  \vmorphposnrt=\z@
\newdimen\vmorphposnlft \vmorphposnlft=\z@

\newdimen\hgrid \hgrid=15pt
\newdimen\vgrid \vgrid=15pt

\newdimen\hchannel  \hchannel=0pt
\newdimen\vchannel  \vchannel=0pt
\newdimen\channelwidth \channelwidth=3pt

\dimendef\@hchannel 0         
\dimendef\@vchannel 2         

\catcode`& = \@oldandcatcode \catcode`@ = \@oldatcatcode

\begin{document}

\title{\large{\textbf{ON THE NORM CONVERGENCE OF NONCONVENTIONAL ERGODIC AVERAGES}}}
\author{Tim Austin}
\date{}

\maketitle


\newenvironment{nmath}{\begin{center}\begin{math}}{\end{math}\end{center}}

\newtheorem{thm}{Theorem}[section]
\newtheorem*{thm*}{Theorem}
\newtheorem{lem}[thm]{Lemma}
\newtheorem{prop}[thm]{Proposition}
\newtheorem{cor}[thm]{Corollary}
\newtheorem*{conj*}{Conjecture}
\newtheorem{dfn}[thm]{Definition}
\newtheorem{ques}[thm]{Question}
\theoremstyle{remark}


\newcommand{\A}{\mathcal{A}}
\newcommand{\B}{\mathcal{B}}
\newcommand{\I}{\mathcal{I}}
\newcommand{\frH}{\mathfrak{H}}
\renewcommand{\Pr}{\mathrm{Pr}}
\newcommand{\s}{\sigma}
\renewcommand{\P}{\mathcal{P}}
\renewcommand{\O}{\Omega}
\renewcommand{\S}{\Sigma}
\newcommand{\T}{\mathrm{T}}
\newcommand{\co}{\mathrm{co}}
\newcommand{\e}{\mathrm{e}}
\newcommand{\eps}{\varepsilon}
\renewcommand{\d}{\mathrm{d}}
\newcommand{\im}{\mathrm{i}}
\renewcommand{\l}{\lambda}
\newcommand{\U}{\mathcal{U}}
\newcommand{\G}{\Gamma}
\newcommand{\g}{\gamma}
\renewcommand{\L}{\Lambda}
\newcommand{\hcf}{\mathrm{hcf}}
\newcommand{\F}{\mathcal{F}}
\renewcommand{\a}{\alpha}
\newcommand{\bbN}{\mathbb{N}}
\newcommand{\bbR}{\mathbb{R}}
\newcommand{\bbZ}{\mathbb{Z}}
\newcommand{\bbQ}{\mathbb{Q}}
\newcommand{\bbT}{\mathbb{T}}
\newcommand{\sfE}{\mathsf{E}}
\newcommand{\id}{\mathrm{id}}
\newcommand{\bb}[1]{\mathbb{#1}}
\newcommand{\fr}[1]{\mathfrak{#1}}
\renewcommand{\bf}[1]{\mathbf{#1}}
\renewcommand{\rm}[1]{\mathrm{#1}}
\renewcommand{\cal}[1]{\mathcal{#1}}
\newcommand{\fin}{\nolinebreak\hspace{\stretch{1}}$\lhd$}

\begin{abstract}
We offer a proof of the following nonconventional ergodic theorem:
\begin{thm*}
If $T_i:\bbZ^r \curvearrowright (X,\S,\mu)$ for $i=1,2,\ldots,d$ are
commuting probability-preserving $\bbZ^r$-actions, $(I_N)_{N\geq 1}$
is a F\o lner sequence of subsets of $\bbZ^r$, $(a_N)_{N\geq 1}$ is
a base-point sequence in $\bbZ^r$ and $f_1,f_2,\ldots,f_d \in
L^\infty(\mu)$ then the nonconventional ergodic averages
\[\frac{1}{|I_N|}\sum_{n \in I_N + a_N} \prod_{i=1}^df_i\circ T_i^n\]
converge to some limit in $L^2(\mu)$ that does not depend on the
choice of $(a_N)_{N \geq 1}$ or $(I_N)_{N\geq 1}$.
\end{thm*}
The leading case of this result, with $r = 1$ and the standard
sequence of averaging sets, was first proved by Tao
in~\cite{Tao08(nonconv)}, following earlier analyses of various more
special cases and related results by Conze and
Lesigne~\cite{ConLes84,ConLes88.1,ConLes88.2}, Furstenberg and
Weiss~\cite{FurWei96}, Zhang~\cite{Zha96}, Host and
Kra~\cite{HosKra01,HosKra05}, Frantzikinakis and Kra~\cite{FraKra05}
and Ziegler~\cite{Zie05}. While Tao's proof rests on a conversion to
a finitary problem, we invoke only techniques from classical ergodic
theory, so giving a new proof of his result.
\end{abstract}

\parskip 7pt

\section{Introduction}

The setting for this work is a collection of $d$ commuting
measure-preserving actions $T_i:\bbZ^r \curvearrowright (X,\S,\mu)$,
$i=1,2,\ldots,d$, on a probability space. We present a proof of the
following result:

\begin{thm}[Convergence of multidimensional nonconventional ergodic averages]\label{thm:nonconvave}
If $T_i:\bbZ^r \curvearrowright (X,\S,\mu)$ for $i=1,2,\ldots,d$ are
commuting probability-preserving $\bbZ^r$-actions, $(I_N)_{N\geq 1}$
is a F\o lner sequence of subsets of $\bbZ^r$, $(a_N)_{N\geq 1}$ is
a base-point sequence in $\bbZ^r$ and $f_1,f_2,\ldots,f_d \in
L^\infty(\mu)$ then the nonconventional ergodic averages
\[\frac{1}{|I_N|}\sum_{n \in I_N + a_N} \prod_{i=1}^df_i\circ T_i^n\]
converge to some limit in $L^2(\mu)$ that does not depend on the
choice of $(a_N)_{N \geq 1}$ or $(I_N)_{N\geq 1}$.
\end{thm}

The case of this result with $r = 1$ and the standard sequence of
averaging sets $I_N + a_N := \{1,2,\ldots,N\}$ was first proved by
Tao in~\cite{Tao08(nonconv)}. Tao proceeds by first demonstrating
the equivalence of this result with a finitary assertion about the
behaviour of the restriction of our functions to large finite pieces
of individual orbits. This, in turn, is easily seen to be equivalent
to a purely finitary result about the behaviour of certain sequences
of averages of $1$-bounded functions on $(\bbZ/N\bbZ)^d$ for very
large $N$, and the bulk of Tao's work then goes into proving this
last result. Interestingly, Towsner has shown in~\cite{Tow07} how
the asymptotic behaviour of these purely finitary averages can be
re-interpreted back into an ergodic-theoretic assertion by building
a suitable `proxy' probability-preserving system from these averages
themselves, using constructions from nonstandard analysis.  Tao's
method of analysis can be extended to the case of individual actions
$T_i$ of a higher-rank $r$ and an arbitrary F\o lner sequence in
$\bbZ^r$, but with the base-point shifts $a_N$ all zero, quite
straightforwardly, but seems to require more work in order to be
extended to a proof for the above base-point-uniform version.

In this paper we shall give a different proof of
Theorem~\ref{thm:nonconvave} that uses only more traditional
infinitary techniques from ergodic theory.  Our method is not
affected by shifting the base points of our averages.  In
particular, we recover a new proof of the base-point-fixed case.

The further special case of Theorem~\ref{thm:nonconvave} in which
$r=1$ and $T_i = T^{a_i}$ for some fixed invertible
probability-preserving transformation $T$ and sequence of integers
$a_1$, $a_2$, \ldots, $a_d$ has been the subject of considerable
recent attention, with complete proofs of this case appearing in
works of Host and Kra~\cite{HosKra05} and of Ziegler~\cite{Zie07}.
These, in turn, build on techniques developed in previous papers for
this or other special cases of the theorem by Conze and
Lesigne~\cite{ConLes84,ConLes88.1,ConLes88.2}, Zhang~\cite{Zha96}
and Host and Kra~\cite{HosKra01}, and also on the analysis by
Furstenberg and Weiss in~\cite{FurWei96} of averages of the form
$\frac{1}{N}\sum_{n=1}^N f\circ T^n\cdot g\circ T^{n^2}$ (which, we
stress, do \emph{not} constitute a special case of
Theorem~\ref{thm:nonconvave} in view of the nonlinearity in the
second exponent).

It is this last paper that first formally introduces the important
notion of `characteristic factors' for a system of averages of
products: in our general setting, these comprise a tuple
$(\Xi_1,\Xi_2,\ldots,\Xi_d)$ of $T$-invariant $\s$-subalgebras of
$\S$ such that, firstly,
\[\frac{1}{|I_N|}\sum_{n \in I_N + a_N} \prod_{i=1}^df_i\circ T_i^n - \frac{1}{|I_N|}\sum_{n \in I_N + a_N} \prod_{i=1}^d\sfE_\mu[f_i\,|\,\Xi_i]\circ T_i^n\to 0\]
in $L^2(\mu)$ as  $N\to\infty$ for any $f_1$, $f_2$, \ldots, $f_d
\in L^\infty(\mu)$ and any choice of $(a_N)_{N\geq 1}$ and $(I_N)_{N
\geq 1}$, so that convergence in general will follow if it can be
established when each $f_i$ is $\Xi_i$-measurable; and secondly such
that these factors have a more precisely-describable structure than
the overall original system, so that the asymptotic behaviour of the
right-hand averages above can be analyzed explicitly.

This proof-scheme has not yet been successfully carried out in the
general setting of the present paper.  The analyses of powers of a
single transformation by Host and Kra and by Ziegler both rely on
achieving a very precise classification of all possible
characteristic factors in the form of `nilsystems', within which
setting a bespoke analysis of the convergence of the relevant
ergodic averages has been carried out separately by
Leibman~\cite{Lei05}. In addition, Frantzikinakis and Kra have shown
in~\cite{FraKra05} that nilsystems re-appear in this r\^ole in the
case of a more general collection of invertible single
transformations $T_i$ under the assumption that each $T_i$ and each
difference $T_iT_j^{-1}$ for $i \neq j$ is ergodic, and they deduce
the restriction of Theorem~\ref{thm:nonconvave} to this case also.
However, without this extra ergodicity hypothesis simple examples
show that any tuple of characteristic factors for our system must be
much more complicated, and no good description of such a tuple is
known.

We note in passing that in the course of their analysis
in~\cite{HosKra05} of the case of powers of a single transformation,
Host and Kra also introduce the following `cuboidal' averages
associated to a single action $S:\bbZ^r\curvearrowright (X,\S,\mu)$:
\[\frac{1}{|I_N|}\sum_{n \in a_N + I_N}\prod_{\eta_1,\eta_2,\ldots,\eta_r \in \{0,1\}}f_\eta\circ S^{\eta_1n_1 + \eta_2n_2 + \ldots + \eta_rn_r}.\]
Using their structural results they are able to prove convergence of
these averages also. This result amounts to a different instance of
our Theorem~\ref{thm:nonconvave}, involving $2^r$ commuting
$\bbZ^r$-actions, by defining $T_\eta^n := S^{\eta_1n_1 + \eta_2n_2+
\ldots + \eta_rn_r}$.

In this paper we shall use the possibility of projecting our input
functions $f_i$ onto special factors only in a rather softer way
than in the works above. Noting that the case $d=1$ of
Theorem~\ref{thm:nonconvave} amounts simply to the von Neumann mean
ergodic theorem, we shall show that, if $d \geq 2$, and under the
assumption that Theorem~\ref{thm:nonconvave} holds for collections
of $d-1$ commuting $\bbZ^r$-actions, then from an arbitrary
$\bbZ^d$-system $(X,\S,\mu,T)$ we can always construct an extension
$(\tilde{X},\tilde{\S},\tilde{\mu},\tilde{T})$ and then a factor
$\tilde{\Xi}$ of that extension such that, interpreting our
nonconventional averages as living inside the larger system
$\tilde{X}$, we may replace the first function $f_1$ with its
projection $\sfE_\mu[f_1\,|\,\tilde{\Xi}]$ in the evaluation of
these averages, and this projection is then of such a form that our
nonconventional averages can be immediately approximated by
nonconventional averages involving only $d-1$ actions.  From this
point a proof of Theorem~\ref{thm:nonconvave} follows quickly by
induction on $d$.

It is interesting to note that this overall scheme of building an
extension to a system with a certain additional property and then
showing that this enables us to project just one of the functions
contributing to our nonconventional averages onto a special factor
of that extension is the same as that followed by Furstenberg and
Weiss in~\cite{FurWei96}.  However, the demands they make on their
extension and the ways in which they then exploit it are very
different from ours, and at the level of finer detail there seems to
be no overlap between the proofs.

In fact, the resulting proof of convergence is much more direct than
those previously discovered for the case of powers of a single
transformation (in addition to avoiding Tao's conversion to a
finitary problem).  This is possibly not so surprising: the
construction we use to build our extended system
$(\tilde{X},\tilde{\S},\tilde{\mu},\tilde{T})$ will typically not
respect any additional algebraic structure among the transformations
$T_i$.  Even if these are powers of a single transformation, in
general the $\tilde{T}_i$ will not be, and thus as far as our proof
is concerned this extra assumption lends us no advantage. This is
symptomatic of an important price that we pay in following our
shorter proof: unlike Host and Kra and Ziegler, we obtain
essentially no additional information about the final form that our
nonconventional averages take.  We suspect that substantial new
machinery will be needed in order to describe these limits with any
precision.

Finally, let us take this opportunity to stress that the
substructures of a system $(X,\S,\mu,T)$ that are responsible for
this complexity in the analysis of nonconventional analysis,
although complicated and difficult to describe, are in a sense very
rare.  This heuristic is made precise in the following observation:
if the action $T$ is chosen generically (using the coarse topology
on the collection of probability-preserving actions on a fixed
Lebesgue space $(X,\S,\mu)$, say), then classical arguments (see,
for example, Chapter 8 of Nadkarni~\cite{Nad98}) show that
generically every $T^\g$ is individually weakly mixing, and in this
case not only can our averages be shown to converge using a rather
shorter argument (due to Bergelson in~\cite{Ber87}), but they
converge simply to the product of the separate averages,
$\prod_{i=1}^d\int_Xf_i\,\d\mu$. We should like to propose a view of
the present paper as a contribution to understanding those rare,
specially structured ways in which the averages associated to our
system can deviate from this `purely random' behaviour.

\textbf{Acknowledgements}\quad My thanks go to Vitaly Bergelson,
John Griesmer, Bernard Host, Bryna Kra, Terence Tao and Tamar
Ziegler for several helpful discussions and communications and to
David Fremlin and an anonymous referee for several constructive
suggestions for improvement.

\section{Some preliminary definitions and results}

Our interest in this paper is with a probability-preserving system
$T:\bbZ^{rd}\curvearrowright (X,\S,\mu)$, for which we we will
always assume that the underlying measurable space is standard
Borel. Inside $\bbZ^{rd}$ we distinguish the subgroups $\G_1 :=
\bbZ^r\times \{0\}^{r(d-1)}$, $\G_2 := \{0\}^r\times \bbZ^r\times
\{0\}^{r(d-2)}$, \ldots and $\G_d := \{0\}^{r(d-1)}\times \bbZ^r$.
Each of these is canonically isomorphic to $\bbZ^r$ when written as
a Cartesian product, as here, and we write $\a_i:\bbZ^r
\stackrel{\cong}{\longrightarrow} \G_i$ for these isomorphisms. We
identify the restrictions $T|_{\G_1}$, $T|_{\G_2}$, \ldots,
$T|_{\G_d}$ with the individual $\bbZ^r$ actions
$T_i^{\a_i(\,\cdot\,)}$, and denote them by $T_1$, $T_2$, \ldots,
$T_d$ respectively. Note that, in this setting of group actions, all
of our transformations are implicitly invertible; routine arguments
easily recover versions of Theorem~\ref{thm:nonconvave} suitable for
collections of commuting non-invertible transformations. We shall
sometimes denote a probability-preserving system alternatively by
$(X,\S,\mu,T)$.

We shall also handle several $\mu$-complete $T$-invariant
$\s$-subalgebras of $\S$. As is a standard in ergodic theory we
shall use the term \textbf{factor} either for such a $\s$-subalgebra
or for a probability-preserving intertwining map $\phi:(X,\S,\mu,T)
\to (Y,\Xi,\nu,S)$; to any such $\phi$ we can associate the
invariant $\s$-subalgebra given by the $\mu$-completion of
$\phi^{-1}[\Xi]$ inside $\S$. Henceforth we shall abusively write
$\phi^{-1}[\Xi]$ for this completed $\s$-algebra.

In particular, within our system we can identify the
\textbf{invariant factor} comprising all $A \in \S$ such that
$\mu(T(A) \triangle A) = 0$. This naturally inherits a
$\bbZ^{rd}$-action from the original system.  We shall denote it by
$\S^T$. More generally, if $\G$ is a subgroup of $\bbZ^{rd}$, we can
identify the factor left invariant by $\{T^\g:\ \g \in \G\}$:
extending the above notation, we shall call this the
\textbf{$T|_\G$-isotropy factor} and write it $\S^{T|_\G}$. We shall
frequently refer to this factor in case $\G$ is the subgroup
$\{\a_i(n) - \a_j(n):\ n\in\bbZ^r\}$ for some $i \neq j$, in which
case we write $\S^{T_i = T_j}$ in place of $\S^{T|_{\rm{im}(\a_i -
\a_j)}}$. It will be centrally important throughout this paper that
if $\G$ is Abelian then the isotropy factors $\S^{T|_\G}$ are
$\bbZ^d$-invariant for all $\G \leq \bbZ^d$; for more general group
actions this invariance holds only if $\G$ is a normal subgroup.

We will assume familiarity with the product measurable space
$(X_1\times X_2\times \cdots \times X_d,\S_1\otimes \S_2\otimes
\cdots \otimes \S_d)$ associated to a family of measurable spaces
$(X_i,\S_i)$, $i=1,2,\ldots,d$. Given measurable maps $\psi_i:X_i
\to Y_i$ between such spaces we shall write $\psi_1\times
\psi_2\times \cdots \times \psi_d$ for their coordinate-wise
product:
\[\psi_1\times \psi_2\times \cdots \times \psi_d(x_1,x_2,\ldots,x_d) :=
(\psi_1(x_1),\psi_2(x_2),\ldots,\psi_d(x_d)).\] More generally, if
$T_i:\bbZ^r\curvearrowright (X_i,\S_i)$ is an action for $i =
1,2,\ldots,d$ then we shall write $T_1\times T_2\times \cdots\times
T_d$ for the action $\bbZ^r\curvearrowright (X_1\times X_2\times
\cdots \times X_d,\S_1\otimes \S_2\otimes \cdots \otimes \S_d)$
given by
\[(T_1\times T_2\times \cdots\times T_d)^n := T^n_1\times
T^n_2\times \cdots\times T^n_d.\] If all the $X_i$ are equal to $X$,
all the $Y_i$ to $Y$ and all the $\psi_i$ to $\psi$ then we shall
abbreviate $\psi \times \psi \times \cdots \times \psi$ to
$\psi^{\times d}$, and similarly for actions.

The construction that we later use for our proof of
Theorem~\ref{thm:nonconvave} will also require the standard notion
of an inverse limit of probability-preserving systems; these are
treated, for example, in Examples 6.3 and Proposition 6.4 of
Glasner~\cite{Gla03}.  In addition to the results contained there,
we need the following simple lemmas.

\begin{lem}[Isotropy factors respect inverse limits]\label{lem:isotropyinvlim} Suppose that \[(X,\S,\mu,T) =
\lim_{m\leftarrow}\,(X^{(m)},\S^{(m)},\mu^{(m)},T^{(m)})\] is an
inverse limit of an increasing sequence of $\bbZ^{rd}$-systems with
connecting maps $\theta^{(m')}_{(m)}:X^{(m')}\to X^{(m)}$ for $m'
\geq m$ and overall projections $\theta_{(m)}:X\to X^{(m)}$, and
that $\G \leq \bbZ^{rd}$. Then
\[\S^{T|_\G} = \bigvee_{m\geq 1}\theta_{(m)}^{-1}[(\S^{(m)})^{T^{(m)}|_\G}].\]
\end{lem}

\textbf{Proof}\quad It is clear that $\S^{T|_\G}\supseteq
\theta_{(m)}^{-1}[(\S^{(m)})^{T^{(m)}|_\G}]$ for every $m \geq 1$,
and therefore that $\S^{T|_\G} \supseteq \bigvee_{m\geq
1}\theta_{(m)}^{-1}[(\S^{(m)})^{T^{(m)}|_\G}]$; it remains to prove
the reverse inclusion.  Thus, suppose that $A \in \S$ is
$T|_\G$-invariant.  Then, by the construction of the inverse limit,
for any $\eps > 0$ we can pick some $m_\eps \geq 1$ and some $A_\eps
\in \theta_{(m_\eps)}^{-1}[\S^{(m_\eps)}]$ with $\mu(A\triangle
A_\eps) < \eps$. This last inequality is equivalent to $\|1_A -
1_{A_\eps}\|_1 < \eps$. Since $A$ is $T|_\G$-invariant it follows
that $\|1_A - 1_{A_\eps}\circ T^\g\|_1 < \eps$ for every $\g \in
\G$; hence, letting $f$ be the ergodic average of $1_{A_\eps}$ under
the action of $T|_\G$, we deduce that $f \in
L^\infty(\mu|_{\theta_{(m)}^{-1}[\S^{(m)}]})$, $f$ is
$T|_\G$-invariant and $\|1_A - f\|_1 < \eps$.  Now taking a
level-set decomposition of $f$ yields $T|_\G$-invariant sets in
$\theta_{(m)}^{-1}[\S^{(m)}]$ that approximate $A$ to within $\eps$.
Since $\eps$ was arbitrary this shows that $A$ lies in
$\bigvee_{m\geq 1}\theta_{(m)}^{-1}[(\S^{(m)})^{T^{(m)}|_\G}]$, as
required. \qed

\begin{lem}[Joins respect inverse limits]\label{lem:meetinvlim}
Suppose that $(X,\S,\mu)$ is a probability space and that for each
$i = 1,2,\ldots,k$ we have a tower of $\s$-subalgebras $\Xi_i^{(0)}
\subseteq \Xi_i^{(1)} \subseteq \ldots \subseteq \S$. Then
\[\bigvee_{m\geq 1}\,(\Xi_1^{(m)}\vee \Xi_2^{(m)}\vee \cdots \vee \Xi_k^{(m)}) = \Big(\bigvee_{m\geq 1}\Xi_1^{(m)}\Big)\vee\Big(\bigvee_{m\geq 1}\Xi_2^{(m)}\Big)\vee\cdots\vee\Big(\bigvee_{m\geq 1}\Xi_k^{(m)}\Big).\]
\end{lem}

\textbf{Proof}\quad For every $m \geq 1$ we have
\begin{multline*}
\Xi_1^{(m)}\vee \Xi_2^{(m)}\vee \cdots \vee \Xi_k^{(m)} \subseteq
\Big(\bigvee_{m\geq 1}\Xi_1^{(m)}\Big)\vee\Big(\bigvee_{m\geq
1}\Xi_2^{(m)}\Big)\vee\cdots\vee\Big(\bigvee_{m\geq
1}\Xi_k^{(m)}\Big)\\ \subseteq \bigvee_{m\geq 1}\,(\Xi_1^{(m)}\vee
\Xi_2^{(m)}\vee \cdots \vee \Xi_k^{(m)})
\end{multline*}
and so taking the limit of the left-hand side above gives the
result. \qed

\section{The Furstenberg self-joining}\label{sec:Fberg}

Central to many of the older ergodic-theoretic analyses of special
cases of Theorem~\ref{thm:nonconvave} is a certain multiple
self-joining of the input $\bbZ^{rd}$-system $(X,\S,\mu,T)$.  Given
such a system and also a F\o lner sequence $(I_N)_{N\geq 1}$ and a
base-point sequence $(a_N)_{N\geq 1}$ we can consider the averages
\[\frac{1}{|I_N|}\sum_{n \in a_N + I_N}\int_X \prod_{i=1}^df_i\circ T_i^n\,\d\mu = \frac{1}{|I_N|}\sum_{n \in a_N + I_N}\int_X f_1\cdot\prod_{i=2}^df_i\circ (T_iT_1^{-1})^n\,\d\mu,\]
and now in view of the right-hand expression above, if we know only
the rank-$(d-1)$ case of Theorem~\ref{thm:nonconvave} then we can
deduce that these averages converge, and it is routine to show
(using the standard Borel nature of $(X,\S)$) that the resulting
limit values define a probability measure $\mu^{\ast d}$ on the
product measurable space $(X^d,\S^{\otimes d})$ by the condition
that
\[\mu^{\ast d}(A_1\times A_2\times \ldots \times A_d) := \lim_{N \to \infty}\frac{1}{|I_N|}\sum_{n \in a_N + I_N}\int_X \prod_{i=1}^d1_{A_i}\circ T_i^n\,\d\mu,\]
where we know that this is independent of the choice of
$(a_N)_{N\geq 1}$ and $(I_N)_{N\geq 1}$. It is now also clear that
this measure $\mu^{\ast d}$ is invariant under the $\bbZ^r$-actions
$S_i := T_i^{\times d}$ for $i = 1,2,\ldots, d$ and also under
$S_{d+1} := T_1 \times T_2\times \ldots \times T_d$. We refer to
$(X^d,\S^{\otimes d},\mu^{\ast d})$ as the \textbf{Furstenberg
self-joining} of the space $(X,\S,\mu)$ associated to the action
$T$, in light of its historical genesis in Furstenberg's work on the
ergodic theoretic approach to Szemer\'edi's Theorem~(\cite{Fur77});
note, in particular, that the one-dimensional marginals of
$\mu^{\ast d}$ on $(X,\S)$ all coincide with $\mu$. Given this
self-joining, we shall write $\pi_1$, $\pi_2$, \ldots, $\pi_d$ for
the projection maps onto the $d$ copies of $(X,\S,\mu)$ that are its
coordinate factors.

In the sequel we will need to work simultaneously with the
Furstenberg self-joinings of a system $(X,\S,\mu,T)$ and an
extension $\psi:(\tilde{X},\tilde{\S},\tilde{\mu},\tilde{T}) \to
(X,\S,\mu,T)$ of that system, in which case we can compute easily
that the map $\psi^{\times d}$ identifies
$(\tilde{X}^d,\tilde{\S}^{\otimes d},\tilde{\mu}^{\ast d})$ as an
extension of $(X^d,\S^{\otimes d},\mu^{\ast d})$, and we shall write
$\tilde{\pi}_1$, $\tilde{\pi}_2$, \ldots, $\tilde{\pi}_d$ for the
coordinate-projections of this larger self-joining.

\section{The proof of nonconventional average convergence}

We prove Theorem~\ref{thm:nonconvave} by induction on $d$.  As
remarked above, the case $d= 1$ is simply the von Neumann mean
ergodic theorem, so let us suppose that $d \geq 2$ and that the
result is known to be true for all systems of at most $d-1$
commuting $\bbZ^r$-actions.

\subsection{Characteristic factors and pleasant systems}

As indicated in the introduction, we shall use a rather simple
instance of the notion of `characteristic factors':

\begin{dfn}[Characteristic factors]\label{dfn:charfactors}
Given a system $T:\bbZ^{rd}\curvearrowright (X,\S,\mu)$, a sequence
of \textbf{characteristic factors} for the nonconventional ergodic
averages associated to $T_1$, $T_2$, \ldots $T_d$ is a tuple
$(\Xi_1,\Xi_2,\ldots,\Xi_d)$ of $T$-invariant $\s$-subalgebras of
$\S$ such that
\[\frac{1}{|I_N|}\sum_{n \in a_N + I_N}\prod_{i=1}^df_i\circ T_i^n - \frac{1}{|I_N|}\sum_{n \in a_N + I_N}\prod_{i=1}^d\sfE_\mu[f_i\,|\,\Xi_i]\circ T_i^n\]
in $L^2(\mu)$ as $N \to \infty$ for any $f_1,f_2,\ldots,f_d \in
L^\infty(\mu)$, F\o lner sequence $(I_N)_{N\geq 1}$ and base-point
sequence $(a_N)_{N\geq 1}$.
\end{dfn}

Many previous results on special cases of
Theorem~\ref{thm:nonconvave} have relied on the identification of a
tuple of characteristic factors that could then be described quite
precisely, in the sense that they can be defined by factor maps of
the original system to certain concrete model systems in which a
more detailed analysis of nonconventional averages is feasible. Most
strikingly, the analysis of Host and Kra in~\cite{HosKra05} and
Ziegler in~\cite{Zie07} show that for powers of a single ergodic
transformation there is a single minimal characteristic factor
(equal to all of the $\Xi_i$ above) that may be identified with a
model given by a $d$-step nilsystem, wherein the convergence of the
nonconventional averages and the form of their limits can be
analyzed in great detail.

Here we shall not be so ambitious.  Various examples show that for a
sufficiently complicated system those functions measurable with
respect to either $\S^{T_1}$ or $\S^{T_i = T_1}$ for some $i =
2,3,\ldots,d$ will behave differently (and, in particular,
contribute nontrivially) should they appear as $f_1$ in our
averages, and so we expect any tuple of characteristic factors to
have $\Xi_1 \supseteq \S^{T_1}\vee\bigvee_{i=2}^d\S^{T_i = T_1}$. In
order to explain our approach, let us first suppose that we are
given a system in which we may actually take this to be our first
characteristic factor, and may simply take $\Xi_i := \S$ for
$i=2,3,\ldots,d$.

\begin{dfn}[Pleasant system]
We shall term a system $(X,\S,\mu,T)$ \textbf{pleasant} if
\[\Big(\S^{T_1}\vee\bigvee_{i=2}^d\S^{T_i = T_1},\S,\S,\ldots,\S\Big)\]
is a tuple of characteristic factors.
\end{dfn}

\textbf{Remark}\quad The idea of conditioning just one of the
functions $f_i$ in our averages onto a nontrivial factor already
appears in Furstenberg and Weiss~\cite{FurWei96}, in whose
terminology such a factor is `partially characteristic'. \fin

The main observation of this subsection is that, given convergence
of nonconventional averages in general for systems of $d-1$ actions,
we can easily deduce that convergence for pleasant systems of $d$
actions. Let us first record separately an elementary robustness
result for nonconventional averages that we shall need shortly.

\begin{lem}\label{lem:contractive}
For any $f_1$, $f_2$, \ldots, $f_d \in L^\infty(\mu)$ and $N \geq 1$
we have
\[\Big\|\frac{1}{|I_N|}\sum_{n \in a_N + I_N}\prod_{i=1}^df_i\circ T_i^n\Big\|_2 \leq \|f_1\|_2\cdot \prod_{i=2}^d\|f_i\|_\infty.\]
\end{lem}

\textbf{Proof}\quad This is clear from the termwise estimate
\[\Big\|\prod_{i=1}^df_i\circ T_i^n\Big\|_2 \leq \|f_1\circ T_1^n\|_2\cdot \prod_{i=2}^d\|f_i\circ T_i^n\|_\infty = \|f_1\|_2\cdot \prod_{i=2}^d\|f_i\|_\infty.\]
and the triangle inequality. \qed

\begin{cor}\label{cor:robustness}
The nonconventional averages
\[\frac{1}{|I_N|}\sum_{n \in a_N + I_N}\prod_{i=1}^df_i\circ T_i^n\]
converge in $L^2(\mu)$ for the $d$-tuple of functions $f_1$, $f_2$,
\ldots, $f_d \in L^\infty(\mu)$ if the corresponding averages are
known to converge for all the $d$-tuples $f^{(m)}_1$, $f_2$, \ldots,
$f_d$ for some sequence $f_1^{(m)} \in L^\infty(\mu)$ that converges
to $f_1$ in $L^2(\mu)$. \qed
\end{cor}

\begin{prop}[Nonconventional average convergence for pleasant
systems]\label{prop:pleasantnonconvave} If
$T:\bbZ^{rd}\curvearrowright (X,\S,\mu)$ is pleasant and
Theorem~\ref{thm:nonconvave} is known to hold for all systems of
$d-1$ commuting actions, then its conclusion also holds for
$(X,\S,\mu,T)$.
\end{prop}

\textbf{Proof}\quad Writing $\Xi :=
\S^{T_1}\vee\bigvee_{i=2}^d\S^{T_i = T_1}$,
Definition~\ref{dfn:charfactors} tells us that
\[\frac{1}{|I_N|}\sum_{n \in a_N + I_N}\prod_{i=1}^df_i\circ T_i^n - \frac{1}{|I_N|}\sum_{n \in a_N + I_N}(\sfE_\mu[f_1\,|\,\Xi]\circ T_1^n)\cdot \prod_{i=2}^df_i\circ T_i^n \to 0\]
for any $f_1$, $f_2$, \ldots, $f_d \in L^\infty(\mu)$, and so it
suffices to prove the desired convergence under the additional
assumption that $f_1$ is $\Xi$-measurable.  However, in this case we
know that we can approximate $f_1$ in $L^2(\mu)$ by finite sums of
the form $\sum_{k=1}^K g_{1,k}\cdot g_{2,k}\cdot\cdots\cdot g_{d,k}$
where $g_{1,k} \in L^\infty(\mu|_{\S^{T_1}})$ and $g_{i,k} \in
L^\infty(\mu|_{\S^{T_1 = T_i}})$ for $i = 2,3,\ldots,d$. Hence by
linearity and Corollary~\ref{cor:robustness} it suffices to prove
convergence for the averages obtained when $f_1$ is replaced by a
single such product:
\[\frac{1}{|I_N|}\sum_{n \in a_N + I_N}((g_1\cdot g_2\cdot\cdots\cdot g_d)\circ T_1^n)\cdot\prod_{i=2}^df_i\circ T_i^n;\]
but now the different invariances that we are assuming for each
$g_i$ imply that $g_1\circ T_1^n = g_1$ and $g_i\circ T_1^n =
g_i\circ T_i^n$ for $i = 2,3,\ldots,d$ and all $n \in \bbZ^r$, and
so the above is simply equal to
\[g_1\cdot \frac{1}{|I_N|}\sum_{n \in a_N + I_N}\prod_{i=2}^d(g_i\cdot f_i)\circ T_i^n.\]
This is a product by the fixed bounded function $g_1$ of a
nonconventional ergodic average associated to the $d-1$ commuting
actions $T_2$, $T_3$, \ldots, $T_d$, and we already know by
inductive hypothesis that these converge in $L^2(\mu)$. This
completes the proof. \qed

Unsurprisingly, there are well-known examples of systems that are
unpleasant: for example, the general $d$-step nilsystems that emerge
in the Host-Kra and Ziegler analyses are such. The simplest example
from among these is the following: if $R_\a$ is an irrational
rotation on $(X,\S,\mu) := (\bbT,\rm{Borel},\rm{Haar})$ and we set
$T_1 := R_\a$, $T_2 := R_{2\a} = T_1^2$, then we can check easily
that $\S^{T_1} = \S^{T_2} = \S^{T_1 = T_2}$ are all trivial, but on
the other hand if $f_2 \in \widehat{\bbT}\setminus\{1_\bbT\}$ and
$f_1 := \overline{f_2}^2$ then $f_1$ and $f_2$ are both orthogonal
to the trivial factor but give
\[\frac{1}{N}\sum_{n=1}^N f_1(T_1^n(t))f_2(T_2^n(t)) = \frac{1}{N}\sum_{n=1}^N \overline{f_2(t)}^2f_2(t)\cdot \overline{f_2(\a)}^2f_2(2\a) \equiv \overline{f_2(t)} \not\to 0\]
as $N\to\infty$.

However, it turns out that we can repair this situation by passing
to a suitable extension.

\begin{prop}[All systems have pleasant extensions]\label{prop:goodextension}
Any $\bbZ^{rd}$-system $(X,\S,\mu,T)$ admits a pleasant extension
$\psi:(\tilde{X},\tilde{\S},\tilde{\mu},\tilde{T}) \to
(X,\S,\mu,T)$.
\end{prop}

From this point, Theorem~\ref{thm:nonconvave} follows at once, since
it is clear that the theorem holds for any system if it holds for
some extension of that system. Proposition~\ref{prop:goodextension}
forms the technical heart of this paper, and we shall prove it in
the next subsection.

\subsection{Building a pleasant extension}\label{subs:pleasant}

We shall build our pleasant extension using the machinery of
Furstenberg self-joinings.  By the remarks of
Section~\ref{sec:Fberg}, given the conclusions of
Theorem~\ref{thm:nonconvave} for systems of $d-1$ commuting
$\bbZ^r$-actions and a system $T:\bbZ^{rd} \curvearrowright
(X,\S,\mu)$ we may form the Furstenberg self-joining
$(X^d,\S^{\otimes d},\mu^{\ast d})$.  Our deduction of pleasantness
for our constructed extension will rest on the following key
estimate.

\begin{lem}[The Furstenberg self-joining controls nonconventional averages]\label{lem:vdC}
If $f_1 \in L^\infty(\mu)$ is such that
\[\int_{X^d} f_1\circ\pi_1\cdot\Big(\prod_{i=2}^df_i\circ \pi_i\Big)\cdot g\,\d\mu^{\ast d} = 0\]
for every choice of $f_2,f_3,\ldots,f_d \in L^\infty(\mu)$ and of
another function $g \in L^\infty(\mu^{\ast d}|_{(\S^{\otimes
d})^{S_{d+1}}})$, then also
\[\frac{1}{|I_N|}\sum_{n \in a_N + I_N}\prod_{i=1}^d f_i\circ T_i^n \to 0\]
in $L^2(\mu)$ for every choice of $f_2,f_3,\ldots,f_d \in
L^\infty(\mu)$ and any F\o lner sequence $(I_N)_{N\geq 1}$ and
base-point sequence $(a_N)_{N\geq 1}$.
\end{lem}

\textbf{Remark}\quad Versions of this result have appeared
repeatedly in previous analyses of more special cases of our main
result; consider, for example, Proposition 5.3 of Zhang~\cite{Zha96}
or Subsection 6.3 of Ziegler~\cite{Zie07}.  The standard proof
applies essentially unchanged in the general setting, and we include
the details here largely for completeness. \fin

\textbf{Proof}\quad Suppose that $f_1$, $f_2$, \ldots, $f_d \in
L^\infty(\mu)$ satisfy the assumptions of the theorem. By the
classical higher-rank van der Corput Lemma (see, for example, the
discussion in Bergelson, McCutcheon and Zhang~\cite{BerMcCZha97})
applied to the bounded $\bbZ^r$-indexed family
$\prod_{i=1}^df_i\circ T_i^n$ in $L^2(\mu)$ we need only prove that
\begin{multline*}\frac{1}{M^{2r}}\sum_{m_1,m_2 \in \{1,2,\ldots,M\}^r}\frac{1}{|I_N|}\sum_{n \in a_N + I_N}\int_X
\prod_{i=1}^d(f_i\circ T_i^{m_1+n} \cdot f_i\circ T_i^{m_2+n})\,\d\mu\\
= \frac{1}{M^{2r}}\sum_{m_1,m_2 \in
\{1,2,\ldots,M\}^r}\frac{1}{|I_N|}\sum_{n \in a_N + I_N}\int_X
\prod_{i=1}^d(f_i\circ T_i^{m_1}\cdot f_i\circ T_i^{m_2})\circ
T_i^n\,\d\mu \to 0
\end{multline*}
as $N \to \infty$ and then $M \to \infty$. However, by the
definition of the Furstenberg self-joining we know that
\begin{multline*}
\frac{1}{|I_N|}\sum_{n \in a_N + I_N}\int_X \prod_{i=1}^d(f_i\circ
T_i^{m_1}\cdot f_i\circ T_i^{m_2})\circ T_i^n\,\d\mu\\ \to \int_X
\prod_{i=1}^d(f_i\cdot f_i\circ
T_i^{m_2-m_1})\circ\pi_i\,\d\mu^{\ast d}
\end{multline*}
as $N \to \infty$.  Now, when we the averages these limiting values
over $m_1$ and $m_2 \in \{1,2,\ldots,M\}^r$, we clearly obtain
convex combinations of uniform averages over increasingly large
ranges of $m_2-m_1$ of the last expression above, and so appealing
to the usual mean ergodic theorem for the $\bbZ^r$-action $S_{d+1}:=
T_1\times T_2\times\cdots \times T_d$ in $L^2(\mu^{\ast d})$ we
deduce that our above double averages converge to
\[\int_X \prod_{i=1}^df_i\circ\pi_i\cdot\Big(\lim_{M\to \infty}\frac{1}{M^r}\sum_{m\in\{1,2,\ldots,M\}^r}\Big(\prod_{i=1}^d f_i\circ \pi_i\Big)\circ S_{d+1}^m\Big)\,\d\mu^{\ast d}.\]
Setting
\[g:= \lim_{M\to \infty}\frac{1}{M^r}\sum_{m \in \{1,2,\ldots,M\}^r}\Big(\prod_{i=1}^d f_i\circ \pi_i\Big)\circ S_{d+1}^m\]
this is precisely an integral of the form that we are assuming
vanishes, as required. \qed

We are now in a position to construct our pleasant extension.

\textbf{Proof of Proposition~\ref{prop:goodextension}}\quad We need
to find an extension $(\tilde{X},\tilde{\S},\tilde{\mu},\tilde{T})$
such that, setting
\[\Xi:= \tilde{\S}^{\tilde{T}_1}\vee \tilde{\S}^{\tilde{T}_2 = \tilde{T}_1}\vee \cdots \vee \tilde{\S}^{\tilde{T}_d = \tilde{T}_1},\]
we have
\[\frac{1}{|I_N|}\sum_{n \in a_N + I_N}\big(\tilde{f}_1 - \sfE_{\tilde{\mu}}[\tilde{f}_1\,|\,\Xi]\big)\circ \tilde{T}_1^n\cdot\prod_{i=2}^d \tilde{f}_i\circ \tilde{T}_i^n \to 0\quad\quad\hbox{in}\ L^2(\tilde{\mu})\]
for any $\tilde{f}_1$, $\tilde{f}_2$, \ldots, $\tilde{f}_d \in
L^\infty(\tilde{\mu})$. By Lemma~\ref{lem:vdC}, this will follow if
we can guarantee instead that
\[\int_{\tilde{X}^d} \tilde{f}_1\circ\tilde{\pi}_1\cdot\Big(\prod_{i=2}^d\tilde{f}_i\circ \tilde{\pi}_i\Big)\cdot \tilde{g}\,\d\tilde{\mu}^{\ast d} = \int_{\tilde{X}^d} \sfE_{\tilde{\mu}}[\tilde{f}_1\,|\,\Xi]\circ\tilde{\pi}_1\cdot\Big(\prod_{i=2}^d\tilde{f}_i\circ \tilde{\pi}_i\Big)\cdot \tilde{g}\,\d\tilde{\mu}^{\ast d}\]
for every choice of $\tilde{f}_1$, $\tilde{f}_2$, \ldots,
$\tilde{f}_d \in L^\infty(\tilde{\mu})$ and of another function
$\tilde{g} \in L^\infty(\tilde{\mu}^{\ast d}|_{(\tilde{\S}^{\otimes
d})^{\tilde{S}_{d+1}}})$.

We shall show that this obtains for the inverse limit of a tower of
extensions of $(X,\S,\mu,T)$ constructed from the Furstenberg
self-joinings themselves.

\textbf{Step 1: construction of the extension}\quad Given the
original system $(X,\S,\mu,T)$ we define an extension
$\psi^{(1)}:(X^{(1)},\S^{(1)},\mu^{(1)},T^{(1)}) \to (X,\S,\mu,T)$
by setting $(X^{(1)},\S^{(1)},\mu^{(1)}) := (X^d,\S^{\otimes
d},\mu^{\ast d})$, $\psi^{(1)} := \pi_1$ and with the
$\bbZ^r$-actions
\begin{eqnarray*}
T^{(1)}_1 &:=& S_{d+1},\\
T^{(1)}_2 &:=& S_2,\\
&\vdots&\\
T^{(1)}_d &:=& S_d
\end{eqnarray*}
(note that we lift $T_1$ to $S_{d+1}$, rather than to $S_1$).  We
may now iterate this construction on the systems that emerge from it
to build a whole tower of extensions
$(X^{(m)},\S^{(m)},\mu^{(m)},T^{(m)}) \to
(X^{(m-1)},\S^{(m-1)},\mu^{(m-1)},T^{(m-1)})$ for $m \geq 1$, where
we set $(X^{(0)},\S^{(0)},\mu^{(0)},T^{(0)}) := (X,\S,\mu,T)$. Note
that since each $(X^{(m+1)},\S^{(m+1)},\mu^{(m+1)})$ is the $d$-fold
Furstenberg self-joining of $(X^{(m)},\S^{(m)},\mu^{(m)})$, in
addition to the factor map $\pi_1^{(m)}$ given by the projection
onto the first coordinate in this self-joining it carries $d-1$
other such maps corresponding to the projections onto the other
coordinates; let us denote these by $\psi_2^{(m)}$, $\psi_3^{(m)}$,
\ldots, $\psi_d^{(m)}$.

We will take $(\tilde{X},\tilde{\S},\tilde{\mu},\tilde{T})$ to be
the inverse limit $\lim_{m
\leftarrow}(X^{(m)},\S^{(m)},\mu^{(m)},T^{(m)})$, and show that this
has the desired property.  Write $\psi:\tilde{X} \to X$ for the
overall factor map back onto the original probability space,
$\theta^{(m')}_{(m)}:X^{(m')} \to X^{(m)}$ for the connecting
projections of our inverse system, and also $\theta_{(m)}:\tilde{X}
\to X^{(m)}$ for the overall projection from the limit system, so
that $\psi = \theta_{(0)}$. Write $\pi^{(m)}_i$ for the coordinate
projections $(X^{(m)})^d \to X^{(m)}$ and $\tilde{\pi}_i$ for the
coordinate projections $\tilde{X}^d \to \tilde{X}$. Finally, let
\[\Xi^{(m)}:= (\S^{(m)})^{T^{(m)}_1}\vee
(\S^{(m)})^{T^{(m)}_1 = T^{(m)}_2}\vee\cdots\vee
(\S^{(m)})^{T^{(m)}_1 =T^{(m)}_d}\] and
\[\Xi:= \tilde{\S}^{\tilde{T}_1}\vee \tilde{\S}^{\tilde{T}_2 = \tilde{T}_1}\vee \cdots \vee \tilde{\S}^{\tilde{T}_d = \tilde{T}_1};\]
combining Lemmas~\ref{lem:isotropyinvlim} and~\ref{lem:meetinvlim}
we deduce that $\Xi = \bigvee_{m \geq 1}
\theta_{(m)}^{-1}[\Xi^{(m)}]$.

We can depict the tower of systems constructed above in the
following commutative diagram:
$$\commdiag{
(\tilde{X},\tilde{\S},\tilde{\mu})& \mapleft^{\tilde{\pi}_1}&
\big(\tilde{X}^d,\tilde{\S}^{\otimes d},\tilde{\mu}^{\ast d}\big)\cr
\mapdown & & \mapdown\cr \vdots &\vdots & \vdots\cr
\mapdown\lft{\theta^{(3)}_{(2)}} & &
\mapdown\lft{(\theta^{(3)}_{(2)})^{\times d}}\cr
(X^{(2)},\S^{(2)},\mu^{(2)})& \mapleft^{\pi^{(2)}_1}&
\big((X^{(2)})^d,(\S^{(2)})^{\otimes d},(\mu^{(2)})^{\ast d}\big)\cr
\mapdown\lft{\theta^{(2)}_{(1)}} & &
\mapdown\lft{(\theta^{(2)}_{(1)})^{\times d}}\cr
(X^{(1)},\S^{(1)},\mu^{(1)})& \mapleft^{\pi^{(1)}_1}&
\big((X^{(1)})^d,(\S^{(1)})^{\otimes d},(\mu^{(1)})^{\ast d}\big)\cr
\mapdown\lft{\theta^{(1)}} & & \mapdown\lft{(\theta^{(1)})^{\times
d}}\cr (X,\S,\mu)& \mapleft^{\pi_1}& (X^d,\S^{\otimes d}\mu^{\ast
d})\cr }$$ where, in addition, by construction we have
\[(X^{(m+1)},\S^{(m+1)},\mu^{(m+1)}) =
((X^{(m)})^d,(\S^{(m)})^{\otimes d},(\mu^{(m)})^{\ast d})\] for very
$m \geq 0$ with the actions $T_i^{(m+1)}$ selected from among the
$S_i^{(m)}$ as above, and under this identification the maps
$\theta^{(m+1)}_{(m)}$ and $\pi^{(m)}_1$ agree. On the other hand,
the maps $\pi^{(m+1)}_1$ and $(\theta^{(m+1)}_{(m)})^{\times d}$ do
not agree.

\textbf{Step 2: proof of pleasantness}\quad We will now prove that
for any $\tilde{f}_1$, $\tilde{f}_2$, \ldots, $\tilde{f}_d \in
L^\infty(\tilde{\mu})$ and $\tilde{g} \in L^\infty(\tilde{\mu}^{\ast
d}|_{(\tilde{\S}^{\otimes d})^{\tilde{S}_{d+1}}})$ we have
\[\int_{\tilde{X}^d} \tilde{f}_1\circ\tilde{\pi}_1\cdot\Big(\prod_{i=2}^d\tilde{f}_i\circ \tilde{\pi}_i\Big)\cdot \tilde{g}\,\d\tilde{\mu}^{\ast d} = \int_{\tilde{X}^d} \sfE_{\tilde{\mu}}[\tilde{f}_1\,|\,\Xi]\circ\tilde{\pi}_1\cdot\Big(\prod_{i=2}^d\tilde{f}_i\circ \tilde{\pi}_i\Big)\cdot \tilde{g}\,\d\tilde{\mu}^{\ast d}\]
By continuity in $L^2(\tilde{\mu})$ and the definition of inverse
limit, we may assume further that there is some finite $m \geq 1$
such that $\tilde{f}_i = f_i\circ\theta_{(m)}$ and $\tilde{g} =
g\circ \theta_{(m)}^{\times d}$ for some $f_1$, $f_2$, \ldots, $f_d
\in L^\infty(\mu^{(m)})$ and $g \in L^\infty((\mu^{(m)})^{\ast
d}|_{(\S^{(m)})^{S^{(m)}_{d+1}}})$. Given this the left-hand
expression above can be re-written at level $m$ as
\[\int_{(X^{(m)})^d} f_1\circ\pi^{(m)}_1\cdot\Big(\prod_{i=2}^df_i\circ \pi^{(m)}_i\Big)\cdot g\,\d(\mu^{(m)})^{\ast d}.\]

For any $m' \geq m$, since $((X^{(m')})^d,(\S^{(m)})^{\otimes
d},(\mu^{(m')})^{\ast d}) =: (X^{(m'+1)},\S^{(m+1)},\mu^{(m'+1)})$,
the left-hand side above can also be re-written as
\begin{eqnarray*}
&&\int_{\tilde{X}^d}(f_1\circ\theta_{(m)}\circ\tilde{\pi}_1)\cdot\Big(\prod_{i=2}^df_i\circ\theta_{(m)}\circ\tilde{\pi}_i\Big)\cdot
(g\circ\theta_{(m)}^{\times d})\,\d\tilde{\mu}^{\ast d}\\
&&=
\int_{(X^{(m')})^d}\big((f_1\circ\theta^{(m')}_{(m)})\circ\pi^{(m')}_1)\cdot\Big(\prod_{i=2}^d(f_i\circ\theta^{(m')}_{(m)})\circ\pi^{(m')}_i\Big)\cdot
(g\circ(\theta^{(m')}_{(m)})^{\times d})\,\d(\mu^{(m')})^{\ast d}\\
&&=
\int_{X^{(m'+1)}}\big((f_1\circ\theta^{(m')}_{(m)})\circ\theta^{(m'+1)}_{(m')})\cdot\Big(\prod_{i=2}^d(f_i\circ\theta^{(m')}_{(m)})\circ\psi^{(m'+1)}_i\Big)\cdot
(g\circ(\theta^{(m')}_{(m)})^{\times d})\,\d\mu^{(m'+1)}\\
&&=
\int_{\tilde{X}}\big((f_1\circ\theta^{(m'+1)}_{(m)})\circ\theta_{(m'+1)}\big)\\
&&\quad\quad\quad\quad\quad\quad
\cdot\Big(\prod_{i=2}^d(f_i\circ\theta^{(m')}_{(m)})\circ\psi^{(m'+1)}_i\circ\theta_{(m'+1)}\Big)\cdot
(g\circ(\theta^{(m')}_{(m)})^{\times
d}\circ\theta_{(m'+1)})\,\d\tilde{\mu}.
\end{eqnarray*}
Now, the function
$(f_i\circ\theta^{(m')}_{(m)})\circ\psi^{(m'+1)}_i$ is invariant
under the $\bbZ^r$-action
\[T_i^{(m')}(T_1^{(m')})^{-1} \times T_i^{(m')}(T_2^{(m')})^{-1} \times \cdots \times \rm{id} \times \cdots \times T_i^{(m')}(T_d^{(m')})^{-1} =: T_i^{(m'+1)}(T_1^{(m'+1)})^{-1}\]
for each $i = 2,3,\ldots,d$, and the function
$g\circ(\theta^{(m')}_{(m)})^{\times d}$ is invariant under
$S_{d+1}^{(m')} =: T_1^{(m'+1)}$, so in the last integral above all
factors save the first are
$\theta_{(m'+1)}^{-1}[\Xi^{(m'+1)}]$-measurable, and so we may
condition $f_1\circ\theta^{(m'+1)}_{(m)}$ onto $\Xi^{(m'+1)}$ and
conclude overall that
\begin{eqnarray*}&&\int_{\tilde{X}^d}(f_1\circ\theta_{(m)}\circ\tilde{\pi}_1)\cdot\Big(\prod_{i=2}^df_i\circ\theta_{(m)}\circ\tilde{\pi}_i\Big)\cdot g\circ\theta_{(m)}^{\times d}\,\d\tilde{\mu}^{\ast d}\\ &&=
\int_{\tilde{X}}\big(\sfE[f_1\circ\theta^{(m'+1)}_{(m)})\,|\,\Xi^{(m'+1)}]\circ\theta_{(m'+1)}\big)\\
&&\quad\quad\quad\quad\quad\quad
\cdot\Big(\prod_{i=2}^d(f_i\circ\theta^{(m')}_{(m)})\circ\psi^{(m'+1)}_i\circ\theta_{(m'+1)}\Big)\cdot
(g\circ(\theta^{(m')}_{(m)})^{\times
d}\circ\theta_{(m'+1)})\,\d\tilde{\mu}.
\end{eqnarray*}
Since
\[\sfE[f_1\circ\theta^{(m')}_{(m)}\,|\,\Xi^{(m')}]\circ\theta_{(m')}
\to \sfE[f_1\circ\theta_{(m)}\,|\,\Xi]\] and hence
\[\sfE[f_1\circ\theta^{(m'+1)}_{(m)}\,|\,\Xi^{(m'+1)}]\circ\theta_{(m'+1)} - \sfE[f_1\circ\theta^{(m')}_{(m)}\,|\,\Xi^{(m')}]\circ\theta_{(m')}\to 0\quad\quad\hbox{in}\ L^2(\tilde{\mu})\ \hbox{as}\ m' \to \infty,\]
we next deduce that
\begin{eqnarray*}&&\int_{\tilde{X}}\big(\sfE[f_1\circ\theta^{(m'+1)}_{(m)})\,|\,\Xi^{(m'+1)}]\circ\theta_{(m'+1)}\big)\\
&&\quad\quad\quad\quad\quad\quad
\cdot\Big(\prod_{i=2}^d(f_i\circ\theta^{(m')}_{(m)})\circ\psi^{(m'+1)}_i\circ\theta_{(m'+1)}\Big)\cdot(g\circ(\theta^{(m')}_{(m)})^{\times
d}\circ\theta_{(m'+1)})\,\d\tilde{\mu}\\
&&- \int_{\tilde{X}}\big(\sfE[f_1\circ\theta^{(m')}_{(m)})\,|\,\Xi^{(m')}]\circ\theta_{(m')}\big)\\
&&\quad\quad\quad\quad\quad\quad
\cdot\Big(\prod_{i=2}^d(f_i\circ\theta^{(m')}_{(m)})\circ\psi^{(m'+1)}_i\circ\theta_{(m'+1)}\Big)\cdot
(g\circ(\theta^{(m')}_{(m)})^{\times
d}\circ\theta_{(m'+1)})\,\d\tilde{\mu}\\ &&\to 0
\end{eqnarray*}
as $m' \to \infty$, and by the law of iterated conditional
expectation this last expression is equal to
\begin{multline*}\int_{\tilde{X}}\big(\sfE\big[\sfE[f_1\circ\theta^{(m')}_{(m)})\,|\,\Xi^{(m')}]\circ\theta^{(m'+1)}_{(m')}\,\big|\,\Xi^{(m'+1)}\big]\circ\theta_{(m'+1)}\big)\\ \cdot\Big(\prod_{i=2}^d(f_i\circ\theta^{(m')}_{(m)})\circ\psi^{(m'+1)}_i\circ\theta_{(m'+1)}\Big)
\cdot (g\circ(\theta^{(m')}_{(m)})^{\times
d}\circ\theta_{(m'+1)})\,\d\tilde{\mu}.
\end{multline*}
However, by exactly analogous reasoning to that above applied with
$m'$ in place of $m$ and the collection of functions $\sfE[f_1\circ
\theta^{(m')}_{(m)}\,|\,\Xi^{(m')}]\circ\theta_{(m')}$,
$f_i\circ\theta_{(m)} =
(f_i\circ\theta^{(m')}_{(m)})\circ\theta_{(m')}$ for $i =
2,3,\ldots,d$ and $g\circ \theta_{(m)}^{\times d} =
(g\circ(\theta^{(m')}_{(m)})^{\times d})\circ\theta_{(m')}^{\times
d}$ we deduce that this is equal to
\begin{multline*}\int_{\tilde{X}^d}(\sfE[f_1\circ\theta^{(m')}_{(m)}\,|\,\Xi^{(m')}]\circ\theta_{(m')}\circ\tilde{\pi}_1)\cdot\Big(\prod_{i=2}^df_i\circ\theta_{(m)}\circ\tilde{\pi}_i\Big)\cdot
g\circ\theta_{(m)}^{\times d}\,\d\tilde{\mu}^{\ast d}\\ \to
\int_{\tilde{X}^d}(\sfE[f_1\circ\theta_{(m)}\,|\,\Xi]\circ\tilde{\pi}_1)\cdot\Big(\prod_{i=2}^df_i\circ\theta_{(m)}\circ\tilde{\pi}_i\Big)\cdot
g\circ\theta_{(m)}^{\times d}\,\d\tilde{\mu}^{\ast d}
\end{multline*}
as $m' \to \infty$, as required. \qed

It is clear that the assertion of Theorem~\ref{thm:nonconvave} must
hold for any system if it holds for some extension of that system,
and so, as remarked previously, it now follows in full generality by
combining Proposition~\ref{prop:pleasantnonconvave} and
Proposition~\ref{prop:goodextension}. \qed

\textbf{Remarks}\quad Intuitively, at each step in our iterative
construction of the tower of extensions \[(X,\S,\mu,T)\leftarrow
(X^{(1)},\S^{(1)},\mu^{(1)},T^{(1)}) \leftarrow
(X^{(2)},\S^{(2)},\mu^{(2)},T^{(2)})\leftarrow \cdots\] we are
introducing a new supply of functions that are invariant under
either $T_1^{(j)}$ or $T_1^{(j)}(T_i^{(j)})^{-1}$ that can
contribute to building a conditional expectation of $f_1$ that will
serve as a good approximation to it for the purpose of evaluating
our integral.  However, at each such step we introduce new functions
on the larger system that we will also then need to handle in this
way, and these will not be taken care of until the next extension.
It is for this reason that the present construction relies on the
passage all the way to an inverse limit.

Considering informally how the pleasant extension enables us to
bring the proof of Proposition~\ref{prop:pleasantnonconvave} to bear
on a more general system, we can locate the concrete appearance of
the extension $(\tilde{X},\tilde{\S},\tilde{\mu},\tilde{T})$ when we
approximate $f_1$ by $\sum_{k=1}^K g_{1,k}\cdot
g_{2,k}\cdot\cdots\cdot g_{d,k}$: the point is that while this sum
overall approximates a function on the smaller system
$(X,\S,\mu,T)$, the individual functions $g_i$ that appear within it
do not, and then when we separately replace composition with
$\tilde{T}_1^n$ by $\tilde{T}_i^n$ for these functions this requires
us to keep track of their individual orbits inside
$L^\infty(\tilde{\mu})$, which will in general not be confined to
$L^\infty(\mu)$. \fin

\section{Discussion}

\subsection{Alternative constructions of the extension}

The scheme we have adopted to construct our pleasant inverse limit
extension $(\tilde{X},\tilde{\S},\tilde{\mu},\tilde{T})$ of
$(X,\S,\mu,T)$ is far from canonical. In particular, there is more
than one way to use some self-joining of $(X,\mu)$ built using the
original transformations $T$ to control the convergence of
nonconventional averages, as we have done with the Furstenberg
self-joining via Lemma~\ref{lem:vdC}. While this choice seems
particularly well-adapted to giving a quick inductive proof of
Theorem~\ref{thm:nonconvave}, it may be instructive to describe
briefly an alternative such self-joining that could be used in a
similar way. This is a simple generalization of the space
$(X^{[d]},\S^{[d]},\mu^{[d]})$ constructed by Host and Kra for their
proof in~\cite{HosKra05} of Theorem~\ref{thm:nonconvave} in the case
of powers of a single transformation.

Given our original system $(X,\S,\mu,T)$, we construct a sequence of
self-joinings $(X^{[1]},\S^{[1]},\mu^{[1]},T^{[1]})$,
$(X^{[2]},\S^{[2]},\mu^{[2]},T^{[2]})$, \ldots,
$(X^{[d]},\S^{[d]},\mu^{[d]},T^{[d]})$, where each
$(X^{[i]},\S^{[i]},\mu^{[i]},T^{[i]})$ is a $2^i$-fold self-joining
of $(X,\S,\mu,T)$, iteratively as follows.  First set
$(X^{[1]},\S^{[1]}) := (X^2,\S^{\otimes 2})$ and let $\mu^{[1]}$ be
the relatively independent self-joining $\mu\otimes_{\S^{T_1}}\mu$
of $\mu$ over the isotropy factor $\S^{T_1}$ (see, for example,
Section 6.1 of Glasner~\cite{Gla03} for the general construction of
relatively independent self-joinings).  In addition, lift $T_1$ to
$T_1\times\rm{id}_{X}$ and $T_i$ to $T^{[1]}_i := T_i\times T_i$ for
$i = 2,3,\ldots,d$. It is clear from our construction that these
preserve $\mu^{[1]}$. Finally, let $\pi_1$ be the projection of
$X^2$ onto the first coordinate. Now to form
$(X^{[2]},\S^{[2]},\mu^{[2]},T^{[2]})$ we apply this construction to
the system $(X^{[1]},\S^{[1]},\mu^{[1]},T^{[1]})$ but taking the
relatively independent self-product of $\mu^{[1]}$ over the
different isotropy factor $\S^{T^{[1]}_1 = T^{[1]}_2}$, and lifting
$T^{[1]}_1$ to $T_1^{[1]}\times T_2^{[2]}$ and $T_i^{[1]}$ to
$T_i^{[1]}\times T_i^{[1]}$ for $i=2,3,\ldots,d$. We continue
iterating this construction, at each step forming
$(X^{[k]},\S^{[k]},\mu^{[k]},T^{[k]})$ by taking the relatively
independent self-product over $\S^{T^{[k-1]}_1 = T^{[k-1]}_i}$ and
lifting $T_1^{[k-1]}$ to $T_1^{[k-1]}\times T_k^{[k-1]}$ and
$T_i^{[k-1]}$ to $T_i^{[k-1]}\times T_i^{[k-1]}$ for
$i=2,3,\ldots,d$, until we reach $k=d$. This gives the Host-Kra
self-joining.  Our convention is to index the $2^d$-fold product
$X^{[d]}$ that results by the power set $\P[d]$ (the set of all
subsets of $\{1,2,\ldots,d\}$), so that $X^{[d]} = X^{\P[d]}$, in
such a way that $X^{[1]}$ corresponds to the factor
$X^{\{\emptyset,\{1\}\}}$ of this larger product, $X^{[2]}$ to the
factor $X^{\{\emptyset,\{1\},\{2\},\{1,2\}\}}$, and so on. In
addition, we write $\pi_\a^{[d]}$ for the $2^d$ coordinate
projections $X^{\P[d]}\to X$. We can now easily concatenate the
above specifications to write out the resulting transformations
$T^{[d]}_i$ in terms of the original $T_j$: $T^{[d]}_1 = \prod_{\a
\in \P[d]}T_{1,\a}$ with
\[T_{1,\a} :=
\left\{\begin{array}{ll}T_1&\quad\hbox{if}\ \a=\emptyset\\\rm{id}_X&\quad\hbox{if}\ \a = \{1\}\\
T_i&\quad\hbox{if}\ \max\a = i\ \hbox{for}\
i=2,3,\ldots,d,\end{array}\right.\] and $T_i^{[d]}$ is simply
$T_i^{\times \cal{P}[d]}$ for $i=2,3,\ldots,d$.

This can serve as an alternative to the Furstenberg self-joining in
light of the following lemma:

\begin{lem}[The Host-Kra self-joining controls nonconventional averages]\label{lem:HKvdC}
If $f_1 \in L^\infty(\mu)$ is such that
\[\int_{X^{\P[d]}} f_1\circ\pi_\emptyset\cdot\Big(\prod_{\a\in\P[d]\setminus\{\emptyset\}}f_\a\circ \pi_\a\Big)\,\d\mu^{[d]} = 0\]
for every choice of $f_\a\in L^\infty(\mu)$ for $\a \in
\P[d]\setminus\{\emptyset\}$, then also
\[\frac{1}{N}\sum_{n=1}^N\prod_{i=1}^d f_i\circ T_i^n \to 0\]
in $L^2(\mu)$ for every choice of $f_2,f_3,\ldots,f_d \in
L^\infty(\mu)$.
\end{lem}

\textbf{Proof}\quad This follows essentially by $d$ times applying
alternately the van der Corput estimate, just as in the proof of
Lemma~\ref{lem:vdC}), and then the Cauchy-Schwarz inequality for the
space $L^2(\mu)$. The argument is just as for the case of powers of
a single transformation treated by Host and Kra in~\cite{HosKra05}
(see their Theorem 12.1 and the construction of Section 4), and we
omit the details. \qed

Writing $(X^{(1)},\S^{(1)},\mu^{(1)},T^{(1)}) :=
(X^{[d]},\S^{[d]},\mu^{[d]},T^{[d]})$ we can now use the machinery
of Host-Kra self-joinings to build a tower of extensions of
$(X,\S,\mu,T)$ and deduce that their inverse limit is pleasant, as
we did using the Furstenberg self-joining in
Proposition~\ref{prop:goodextension}. This requires grouping
together the various factors in the integrand of
\[\int_{X^{\P[d]}} f_1\circ\pi_\emptyset\cdot\Big(\prod_{\a\in\P[d]\setminus\{\emptyset\}}f_\a\circ \pi_\a\Big)\,\d\mu^{[d]} = 0\]
according to the partition $\P[d]\setminus\{\emptyset\} =
\bigcup_{i=1}^d\{\a:\ \max\a = i\}$, noting that the above explicit
description of $T_1^{[d]}$ tells us that $f_{\{1\}}\circ
\pi^{[d]}_{\{1\}}$ is $T_1^{[d]}$-invariant and that
\[\prod_{\a:\ \max\a = i}f_\a\circ\pi^{[d]}_\a\]
is $T_1^{[d]}(T_i^{[d]})^{-1}$-invariant for $i=2,3,\ldots,d$. The
remaining details of the argument are almost identical to those for
Proposition~\ref{prop:goodextension}. We note that in this argument
the one-step extension $(X^{(1)},\S^{(1)},\mu^{(1)},T^{(1)})$ is
already the top member of a height-$d$ tower of self-joinings. These
two towers serve different purposes in the proof, and should not be
confused: the $d$ smaller extensions used to build up to
$(X^{(1)},\S^{(1)},\mu^{(1)},T^{(1)})$ correspond to the $d$ appeals
to the van der Corput estimate during the proof of
Lemma~\ref{lem:HKvdC}.

The choice between the Furstenberg and Host-Kra self-joinings
certainly affects the structure of the pleasant extension that
emerges, but seems to make little difference to the overall
complexity of the proof, since we do not exploit any of this more
particular structure.  The advantage of the Host-Kra self-joining is
that it does not require an iterative appeal to
Theorem~\ref{thm:nonconvave} for its proof, but on the other hand
that is traded off into a more complicated, alternating appeal to
the van der Corput estimate and the Cauchy-Schwarz inequality in the
proof of Lemma~\ref{lem:HKvdC}, rather than the simple single
application made to prove Lemma~\ref{lem:vdC}.

Looking beyond the above considerations, it may be interesting to
search for a quicker way to pass directly to a pleasant extension:

\textbf{Question}\quad Can we construct a pleasant extension
$(\tilde{X},\tilde{\S},\tilde{\mu},\tilde{T})$ in a finite number of
steps, without invoking an inverse limit? \fin

\textbf{Remark}\quad Since a preprint of this paper first appeared,
Bernard Host has shown in~\cite{Hos??} that by using the above
Host-Kra self-joining, one iteration of the above construction
suffices to produce a pleasant system: the passage to the inverse
limit is already superfluous!  His proof of this requires a slightly
more delicate analysis than the work of our
Subsection~\ref{subs:pleasant}, but in fact it seems likely that it
applies equally well to both self-joinings. \fin

\subsection{Possible further questions}

During the course of proving Theorem~\ref{thm:nonconvave} we have
made essential use of the commutativity of $\bbZ^r$, in addition to
the commutativity of the different actions $T_1$, $T_2$, \ldots,
$T_d$. It is possible that our theorem could be generalized by
considering the averages
\[\frac{1}{|I_N|}\sum_{\g \in a_NI_N}\prod_{i=1}^d f_i\circ T_i^\g\]
for $d$ commuting actions $T_1$, $T_2$, \ldots, $T_d$ on
$(X,\S,\mu)$ of a more general amenable group $\G$ with a F\o lner
sequence $(I_N)_{N\geq 1}$ and base-point sequence $(a_N)_{N\geq
1}$. In this case, if we mimic our straightforward construction of
the Furstenberg self-joining, we obtain a measure $\mu^{\ast d}$ on
$X^d$ that is $T_1\times T_2\times \ldots \times T_d$-invariant, but
it may not now be invariant under any of the diagonal actions
$T_i^{\times d}$. It seems that that ideas of the present paper
cannot yield this stronger result (if it is true at all) without
some additional new insight.

Another generalization of Theorem~\ref{thm:nonconvave} has been
conjectured by Bergelson and Leibman in~\cite{BerLei02}:
\begin{conj*}[Nilpotent nonconventional ergodic averages]
If $T:\G \curvearrowright (X,\S,\mu)$ is a probability-preserving
action of a discrete nilpotent group $\G$ and $\g_1,\g_2,\ldots,\g_d
\in \G$ then for any $f_1,f_2,\ldots,f_d \in L^\infty(\mu)$ the
nonconventional ergodic averages
\[\frac{1}{N}\sum_{n=1}^N \prod_{i=1}^df_i\circ T^{\g_i^n}\]
converge to some limit in $L^2(\mu)$.
\end{conj*}

I do not know whether the methods of the present paper can be
brought to bear on this conjecture; it seems likely that
considerable further new machinery would be needed here also.

In a different direction, it is unknown whether
Theorem~\ref{thm:nonconvave} holds with pointwise convergence in
place of convergence in $L^2(\mu)$. The methods of the present paper
seem to contribute very little to our understanding of this problem;
crucially, while the Furstenberg self-joining allows us to prove
that $f_1 - \sfE_\mu[f_1\,|\,\Xi]$ contributes negligibly to the
$L^2(\mu)$ convergence of our averages inside the extended system,
so that we can replace $f_1$ with $\sfE_\mu[f_1\,|\,\Xi]$, we
currently know of no good way to control this approximation
pointwise, as would be essential for any approach to the question of
pointwise convergence using the machinery of pleasant extensions and
their factors.

\parskip 0pt

\bibliographystyle{abbrv}
\bibliography{nonconvergave4}

\parskip 0pt

\vspace{10pt}

\textsc{Department of Mathematics\\ University of California, Los
Angeles,\\ Los Angeles, CA 90095-1555, USA}

\vspace{7pt}

Email: \verb|timaustin@math.ucla.edu|

Web: \verb|http://www.math.ucla.edu/~timaustin|

\vspace{7pt}

\end{document}